\documentclass[11pt,draftclsnofoot,onecolumn]{IEEEtran}

\usepackage[T1]{fontenc}
\usepackage[utf8]{inputenc}
\usepackage{amsmath,amssymb,amsfonts,mathtools}
\usepackage{amsthm}
\usepackage{booktabs}
\usepackage{array}
\usepackage{graphicx}
\usepackage[numbers,sort&compress]{natbib}
\usepackage{microtype}
\usepackage{url}

\newtheorem{definition}{Definition}
\newtheorem{assumption}{Assumption}
\newtheorem{proposition}{Proposition}

\newtheorem{theorem}{Theorem}
\newtheorem{corollary}{Corollary}

\newcommand{\R}{\mathbb{R}}
\newcommand{\E}{\mathbb{E}}

\newcommand{\cint}{C_{\mathrm{int}}}
\newcommand{\creg}{C_{\mathrm{reg}}}

\newcommand{\opnorm}[1]{\left\lVert #1 \right\rVert_2}
\newcommand{\Frob}[1]{\left\lVert #1 \right\rVert_F}

\title{Control-Channel Informativity for Koopman EDMDc under Behavior-Policy Data}

\author{Yue~Wu%
\thanks{Yue Wu is with the School of Automation Science and Engineering, Xi'an Jiaotong University, Xi'an 710049, China, and with Xinjiang Cigarette Factory, Hongyun Honghe Tobacco (Group) Co., Ltd., Urumqi 830000, China.}}

\begin{document}
\maketitle

\begin{abstract}
Extended dynamic mode decomposition with control (EDMDc) is often trained from trajectories generated by a behavior policy or a pre-existing feedback controller. Such data can predict the observed behavior accurately while failing to identify how new input commands change the lifted state. This paper studies that failure as a control-channel informativity problem. We introduce a conditional intervention certificate, defined as the residual input covariance after projecting the input data away from the active lifted-state feature span. The certificate is the Schur complement of the lifted-state block in the EDMDc information matrix. We prove that its strict positivity is necessary and sufficient for finite-sample sample-identifiability of the lifted control-channel block. If the certificate vanishes, distinct lifted models agree on every collected transition but disagree under counterfactual inputs. We then give a closed-loop statistical bound using predictable regressors, conditionally sub-Gaussian transition noise, and a regularized Schur complement. A scalar feedback example shows the unavoidable scaling: under dithered feedback, residual intervention information grows quadratically with the dither amplitude and the control-channel error decreases with the inverse intervention signal-to-noise scale. New experiments verify these scalings exactly in a linear system and diagnostically in controlled Duffing and Van der Pol benchmarks. A larger EDMDc acquisition grid further shows that state coverage, joint regression conditioning, and intervention excitation are complementary diagnostics rather than interchangeable performance scores.
\end{abstract}


\section{Introduction}

Koopman operator methods lift nonlinear dynamics to observable coordinates in which the evolution can be approximated linearly. Extended dynamic mode decomposition and EDMD with control have become standard finite-dimensional tools for prediction and model predictive control \cite{williams2015datadriven,proctor2016dynamic,korda2018linear,brunton2022modern}. Recent work has moved toward control guarantees, including active learning of dynamics \cite{abraham2019active}, generalized excitation conditions \cite{boddupalli2019koopman}, Koopman versions of Willems' fundamental lemma \cite{shang2024willems}, data-driven MPC with EDMD stability margins \cite{bold2025mpc}, feedback design with Koopman models \cite{strasser2025feedback}, and non-asymptotic EDMD or kernel EDMD error analysis \cite{philipp2025error,bold2025kernel}.

The quality of the data used for EDMDc is therefore a control issue, not merely a regression issue. Classical persistent excitation and experiment design quantify whether regressors or inputs are rich enough for identification \cite{narendra1987persistent,willems2005note,hjalmarsson2005experiment,decock2016doptimal,wilson2014trajectory}. Data informativity asks what control properties can be inferred from a finite dataset \cite{van2020data}. Space-filling and active-design methods improve coverage for nonlinear identification \cite{lee2021optimal,kiss2024spacefilling,smits2024spacefilling,liu2025spacefilling,rickenbach2024active}. These tools are necessary, but EDMDc has an additional block-structured question: after the lifted state has explained what it can, is there any independent input variation left to identify the lifted control channel?

The question is concrete. Suppose data are collected under a stabilizing behavior controller. The observed input may vary, and the EDMDc model may predict the closed-loop trajectory well. Nevertheless, if the input is almost a deterministic function of the lifted state, the contribution of the input can be absorbed into the autonomous lifted dynamics on the sampled data. A model trained on such data can fit the behavior trajectory while assigning a wrong response to a new command that breaks the feedback relation. Figure~\ref{fig:concept} summarizes this distinction.

\begin{figure}[t]
\centering
\fbox{\begin{minipage}{0.92\columnwidth}
\centering
\vspace{2mm}
\begin{tabular}{ccc}
behavior data & residualization & intervention test\\
$u_k\approx \mu(\psi(x_k))$ & $U_\perp=U(I-P_Z)$ & new $u$ not on behavior policy\\[1mm]
\multicolumn{3}{c}{$\cint=\lambda_{\min}(N^{-1}U_\perp U_\perp^\top)$}
\end{tabular}
\vspace{2mm}
\end{minipage}}
\caption{EDMDc can fit behavior-policy transitions while leaving the response to new commands unidentified. The proposed certificate measures input variation that remains after conditioning on lifted-state features.}
\label{fig:concept}
\end{figure}

This paper treats the issue as \emph{control-channel informativity}. The central quantity is a conditional intervention certificate $\cint$, the smallest eigenvalue of the input covariance after projection away from active lifted-state features. The object is mathematically a partial-regression Schur complement. The contribution is not the Schur complement alone; it is the control interpretation, the necessary-and-sufficient EDMDc channel identifiability result, the closed-loop self-normalized error bound, and the experiments showing the predicted dither and sample-size scaling.

The contributions are:
\begin{itemize}
\item We define $\cint$ for active, standardized EDMDc data and show that it is the Schur complement controlling the least-squares $B$-block.
\item We prove a finite-sample theorem: the lifted control channel is sample-identifiable if and only if $\cint>0$. If $\cint=0$, there are distinct lifted models that fit all collected transitions but disagree under counterfactual inputs.
\item We upgrade the closed-loop discussion to a formal adapted-noise theorem. Under predictable regressors and conditional sub-Gaussian transition noise, the $B$-block error is governed by a regularized Schur complement.
\item We give a scalar dithered-feedback proposition and experiments showing $\cint\propto \epsilon^2$ and median control-channel error scaling approximately as $(\epsilon\sqrt N)^{-1}$.
\item We validate the diagnostic on nonlinear Duffing and Van der Pol feedback data and retain a broader three-system EDMDc acquisition grid as complementary evidence.
\end{itemize}

The claims are deliberately finite-dimensional and dictionary-dependent. The certificate does not guarantee long-horizon prediction, nonlinear closed-loop optimality, or robustness to dictionary misspecification. It certifies a narrower property: whether the collected EDMDc data contain information that separates the lifted control channel from autonomous lifted evolution.

\section{EDMDc and the Conditional Certificate}

Consider a controlled nonlinear system $x_{k+1}=f(x_k,u_k)$ with $x_k\in\R^n$ and $u_k\in\R^m$. Let $z_k=\psi_A(x_k)\in\R^d$ denote the active lifted dictionary after removing constant or numerically inactive coordinates. The certificates below are computed after centering and standardizing the active lifted features and inputs, using statistics from the identification data.

Given $N$ transitions, define column-sample matrices
\begin{equation}
Z=\begin{bmatrix}z_0&z_1&\cdots&z_{N-1}\end{bmatrix}\in\R^{d\times N},
\end{equation}
\begin{equation}
U=\begin{bmatrix}u_0&u_1&\cdots&u_{N-1}\end{bmatrix}\in\R^{m\times N},
\qquad
Y=\begin{bmatrix}z_1&z_2&\cdots&z_N\end{bmatrix}\in\R^{q\times N}.
\end{equation}
Usually $q=d$. EDMDc fits
\begin{equation}
\min_{A,B}\ \Frob{Y-AZ-BU}^2,
\label{eq:edmdc}
\end{equation}
where $A\in\R^{q\times d}$ and $B\in\R^{q\times m}$. Let $\Phi=[Z^\top,U^\top]^\top$.

\begin{definition}[Joint regression certificate]
The joint EDMDc regression certificate is
\begin{equation}
\creg=\lambda_{\min}\left(\frac{1}{N}\Phi\Phi^\top\right).
\end{equation}
It measures conditioning of the complete predictor $[A,B]$.
\end{definition}

\begin{definition}[Conditional intervention certificate]
Let
\begin{equation}
P_Z=Z^\top(ZZ^\top)^\dagger Z,\qquad M_Z=I_N-P_Z,
\end{equation}
and define $U_\perp=UM_Z$. The conditional intervention certificate is
\begin{equation}
\cint=\lambda_{\min}\left(\frac{1}{N}U_\perp U_\perp^\top\right).
\label{eq:cint}
\end{equation}
\end{definition}

The quantity $\cint$ is not raw input variance. It is the input energy left after removing the component that is linearly predictable from the active lifted state. A small value means that the data do not separate the control channel from the lifted autonomous dynamics. The practical remedy is to alter the experiment, for example by adding dither to the behavior controller, diversifying references, or collecting targeted open-loop interventions over the relevant state region.

\begin{proposition}[Schur complement for the control block]
\label{prop:schur}
For the empirical information matrix
\begin{equation}
G=\frac{1}{N}\Phi\Phi^\top
=\begin{bmatrix}G_{zz}&G_{zu}\\G_{uz}&G_{uu}\end{bmatrix},
\end{equation}
the conditional intervention matrix satisfies
\begin{equation}
\frac{1}{N}U_\perp U_\perp^\top
=G_{uu}-G_{uz}G_{zz}^{\dagger}G_{zu}.
\label{eq:schur}
\end{equation}
If $G_{zz}\succ0$, then $G\succ0$ if and only if $N^{-1}U_\perp U_\perp^\top\succ0$. When $\Phi\Phi^\top\succ0$, the lower-right inverse block of $(\Phi\Phi^\top)^{-1}$ is $(U_\perp U_\perp^\top)^{-1}$.
\end{proposition}

The proof is the standard Schur-complement identity after expanding $UM_ZU^\top$. This is precisely why $\cint$ is the right block-specific certificate: in homoscedastic least squares, the covariance of the fitted $B$-block is controlled by the inverse residual input covariance, not by the raw input covariance.

\section{Deterministic Control-Channel Informativity}

Assume the active lifted data satisfy
\begin{equation}
Y=A_\star Z+B_\star U+R+E,
\label{eq:model}
\end{equation}
where $A_\star$ and $B_\star$ are the finite-dimensional lifted operators in the chosen coordinates, $R$ is a finite-dictionary residual, and $E$ is lifted transition-output noise. We call the lifted control channel sample-identifiable if every pair $(A,B)$ that gives the same noiseless sample predictions has the same $B$-block.

\begin{theorem}[Finite-sample control-channel informativity]
\label{thm:deterministic}
Let $\beta_N=\cint$. The following statements hold.
\begin{enumerate}
\item The lifted control channel $B$ is sample-identifiable if and only if $\beta_N>0$.
\item If $\beta_N>0$, then the $B$-block of every EDMDc least-squares solution is unique and
\begin{equation}
\widehat B=YM_ZU_\perp^\top(U_\perp U_\perp^\top)^{-1}.
\label{eq:bhat}
\end{equation}
Consequently,
\begin{equation}
\widehat B-B_\star=(R+E)M_ZU_\perp^\top(U_\perp U_\perp^\top)^{-1}.
\label{eq:berr}
\end{equation}
\item If $\beta_N=0$, then there exist $a\in\R^m$, $a\ne0$, and $h\in\R^d$ such that $a^\top U=h^\top Z$. For any $c\in\R^q$, define
\begin{equation}
A_c=A_\star-ch^\top,\qquad B_c=B_\star+ca^\top .
\end{equation}
Then $A_cZ+B_cU=A_\star Z+B_\star U$ on the collected data, while for a new pair $(z,u)$,
\begin{equation}
A_cz+B_cu-(A_\star z+B_\star u)=c(a^\top u-h^\top z).
\label{eq:counterfactual}
\end{equation}
\end{enumerate}
\end{theorem}

Theorem~\ref{thm:deterministic} separates two notions that are often conflated. Behavior prediction asks whether $AZ+BU$ matches the collected transitions. Control-channel informativity asks whether the same data determine the map $u\mapsto Bu$ for inputs that do not satisfy the behavior relation. If $\cint=0$, the behavior data cannot answer that intervention question.

Under a fixed-design sub-Gaussian noise model, \eqref{eq:berr} gives the familiar error scale
\begin{equation}
\opnorm{\widehat B-B_\star}
\le
c\sigma\sqrt{\frac{m+q+\log(1/\delta)}{N\beta_N}}
+\frac{\varepsilon_\psi}{\sqrt{\beta_N}},
\label{eq:fixeddesign}
\end{equation}
with probability at least $1-\delta$, provided $\opnorm{RM_Z}\le\sqrt N\varepsilon_\psi$. This bound is useful but incomplete for closed-loop data, because regressors may depend on past disturbances. The next section addresses that point directly.

\section{Closed-Loop Statistical Bound}

Let $\phi_k=[z_k^\top,u_k^\top]^\top\in\R^{d+m}$ and $\Phi=[\phi_0,\ldots,\phi_{N-1}]$. Define the ridge estimator
\begin{equation}
\widehat K_\lambda
=Y\Phi^\top(\Phi\Phi^\top+\lambda I)^{-1},
\qquad
\widehat K_\lambda=[\widehat A_\lambda,\widehat B_\lambda].
\end{equation}
Let $K_\star=[A_\star,B_\star]$ and $V_N=\lambda I+\Phi\Phi^\top$.

\begin{assumption}[Adapted transition noise]
\label{ass:adapted}
Let $\mathcal G_k$ be the sigma-field generated by all variables available before the transition noise $e_k$ is realized. The regressor $\phi_k$ is $\mathcal G_k$-measurable, $\E[e_k\mid\mathcal G_k]=0$, and for every unit vector $v\in\R^q$ and all $\alpha\in\R$,
\begin{equation}
\E[\exp\{\alpha v^\top e_k\}\mid\mathcal G_k]
\le
\exp(\alpha^2\sigma^2/2).
\end{equation}
\end{assumption}

Define the regularized Schur complement
\begin{equation}
S_{u|z,\lambda}
=UU^\top+\lambda I_m
-UZ^\top(ZZ^\top+\lambda I_d)^{-1}ZU^\top .
\label{eq:regschur}
\end{equation}

\begin{theorem}[Closed-loop adapted-noise $B$-block bound]
\label{thm:martingale}
Suppose \eqref{eq:model} holds and Assumption~\ref{ass:adapted} is satisfied. There is an absolute constant $c>0$ such that, with probability at least $1-\delta$,
\begin{equation}
\opnorm{\widehat B_\lambda-B_\star}
\le
\frac{
c\sigma\sqrt{q\log\!\left(\frac{\det(V_N)^{1/2}}{\lambda^{(d+m)/2}\delta}\right)}
+\opnorm{R\Phi^\top V_N^{-1/2}}
+\sqrt{\lambda}\opnorm{K_\star}
}{
\sqrt{\lambda_{\min}(S_{u|z,\lambda})}
}.
\label{eq:martbound}
\end{equation}
\end{theorem}

Theorem~\ref{thm:martingale} makes the closed-loop role of the certificate explicit. Temporal dependence changes the concentration tool from fixed-design matrix concentration to a self-normalized martingale inequality \cite{abbasi2011improved,simchowitz2018learning}. It does not change the block geometry: the $B$-block is still governed by a conditional input Schur complement. The regularization is also useful in finite samples because it avoids singular inverse blocks and exposes the ridge bias term.

\section{Scalar Confounding and Dither Scaling}

The simplest example is the scalar system
\begin{equation}
x_{k+1}=a x_k+b u_k+e_k.
\end{equation}
If the data are collected under deterministic feedback $u_k=\kappa x_k$, then only the closed-loop coefficient $a+b\kappa$ is observed. For any $\widetilde b$, choosing $\widetilde a=a+\kappa(b-\widetilde b)$ gives the same behavior prediction on the collected data. The control coefficient $b$ is not identified.

Now let
\begin{equation}
u_k=\kappa x_k+\epsilon v_k,
\end{equation}
where $v_k$ has zero mean, unit variance, and is independent of $x_k$ and $e_k$. With dictionary $z=x$, the population residual input variance is
\begin{equation}
\E[u_\perp^2]=\epsilon^2.
\end{equation}
Thus the intervention information grows quadratically with dither amplitude.

\begin{proposition}[Unavoidable dither scaling]
\label{prop:scaling}
In the scalar model above, with Gaussian transition noise of variance $\sigma^2$ and dither variance one, the Fisher information for $b$ conditional on $x$ is $N\epsilon^2/\sigma^2$. Consequently, any unbiased estimator satisfies
\begin{equation}
\operatorname{var}(\widehat b)\ge \frac{\sigma^2}{N\epsilon^2}.
\label{eq:crlb}
\end{equation}
When $\epsilon=0$, $b$ is not identifiable from behavior data.
\end{proposition}

Proposition~\ref{prop:scaling} is the statistical counterpart of Theorem~\ref{thm:deterministic}. It shows that the certificate is not merely a renamed conditioning number: small residual intervention excitation produces an unavoidable error amplification for the control channel.

\section{Certificate-to-Feedback Margin}

The main purpose of $\cint$ is identification, not direct nonlinear stability. Still, it can be used in a standard robustness calculation. Assume $q=d$ and consider lifted linear feedback $u=Lz$. Let
\begin{equation}
\widehat F_L=\widehat A+\widehat B L,\qquad
F_{\star,L}=A_\star+B_\star L.
\end{equation}
Let $\eta_A\ge\opnorm{\widehat A-A_\star}$ be obtained from a joint regression bound and $\eta_B\ge\opnorm{\widehat B-B_\star}$ from Theorem~\ref{thm:deterministic} or \ref{thm:martingale}. Define $\eta_{\rm cl}=\eta_A+\opnorm{L}\eta_B$.

\begin{corollary}[Lifted feedback robustness]
\label{cor:feedback}
Suppose there exist $P\succ0$ and $\alpha>0$ such that
\begin{equation}
\widehat F_L^\top P\widehat F_L-P\preceq-\alpha I.
\end{equation}
If
\begin{equation}
2\opnorm{P}\opnorm{\widehat F_L}\eta_{\rm cl}
+\opnorm{P}\eta_{\rm cl}^2
\le \frac{\alpha}{2},
\end{equation}
then $F_{\star,L}^\top P F_{\star,L}-P\preceq-\alpha I/2$.
\end{corollary}

This corollary is a lifted finite-dictionary margin statement. It becomes a nonlinear stability statement only under additional Koopman-invariance or bounded-residual assumptions. Its value here is to show where the $B$-block uncertainty enters a control calculation.

\section{Experiments}

The experiments are designed to answer three questions. First, does the scalar theory predict the observed dither scaling? Second, in nonlinear EDMDc data, can behavior prediction remain accurate while counterfactual response is wrong when $\cint$ is small? Third, how does the new certificate relate to broader sampling diagnostics from the original acquisition grid?

\subsection{Scalar closed-loop scaling}

We simulate
\begin{equation}
x_{k+1}=0.85x_k+0.6u_k+e_k,\qquad e_k\sim\mathcal N(0,0.02^2),
\end{equation}
under $u_k=-0.7x_k+\epsilon v_k$. The sweep uses $\epsilon\in\{0,0.01,0.03,0.1,0.3,1.0\}$, $N\in\{50,100,200,400,800\}$, and 50 seeds. The model is ordinary least squares on $(x_k,u_k)$. Counterfactual error is measured on independent states and uniformly sampled inputs.

\begin{figure*}[t]
\centering
\includegraphics[width=0.98\textwidth]{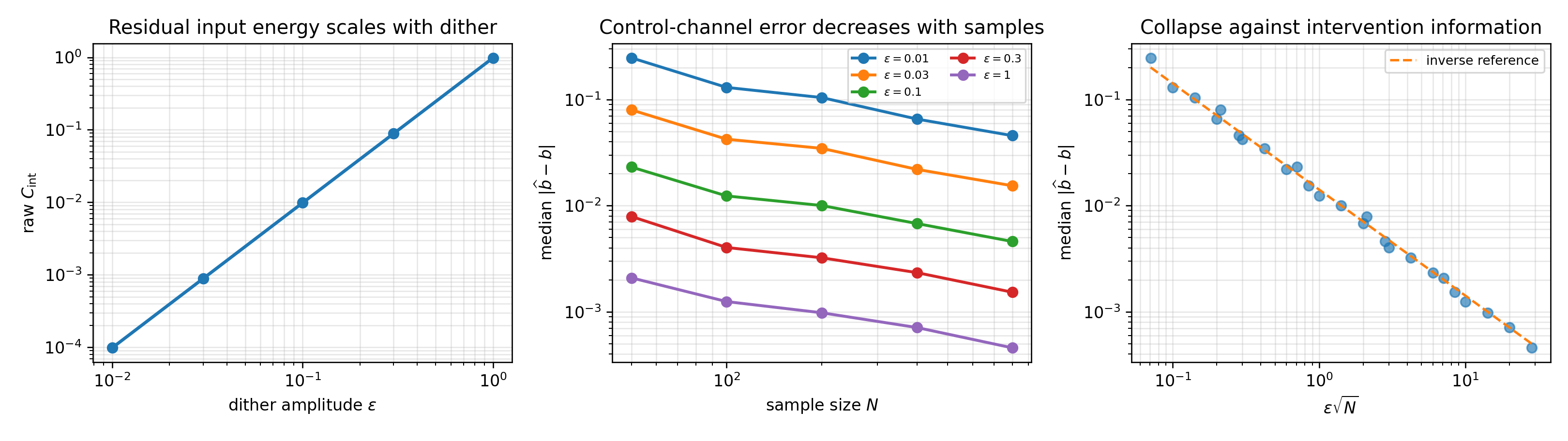}
\caption{Scalar closed-loop scaling experiment. Residual intervention information grows quadratically with dither amplitude, and the control-channel error collapses against $\epsilon\sqrt N$, matching the predicted inverse scaling.}
\label{fig:scalar}
\end{figure*}

Figure~\ref{fig:scalar} directly verifies Proposition~\ref{prop:scaling}. The fitted log-log slope for raw $\cint$ versus $\epsilon$ is $1.9997$, and the fitted slope for median $|\widehat b-b|$ versus $\epsilon\sqrt N$ is $-1.0230$. Deterministic feedback has essentially zero residual intervention information and large counterfactual error despite small behavior error.

\subsection{Nonlinear feedback and dither}

The nonlinear experiments use controlled Duffing and Van der Pol systems with polynomial dictionaries. Data are collected from clipped linear feedback plus Gaussian dither. Budgets are $\{20,40,80\}$ trajectory segments of length 12, dither scales are $\{0,0.02,0.05,0.1,0.2,0.5\}$ times the input bound, and each setting uses 20 seeds. EDMDc uses the same standardized dictionary and ridge selection as the V9 grid.

\begin{figure*}[t]
\centering
\includegraphics[width=0.88\textwidth]{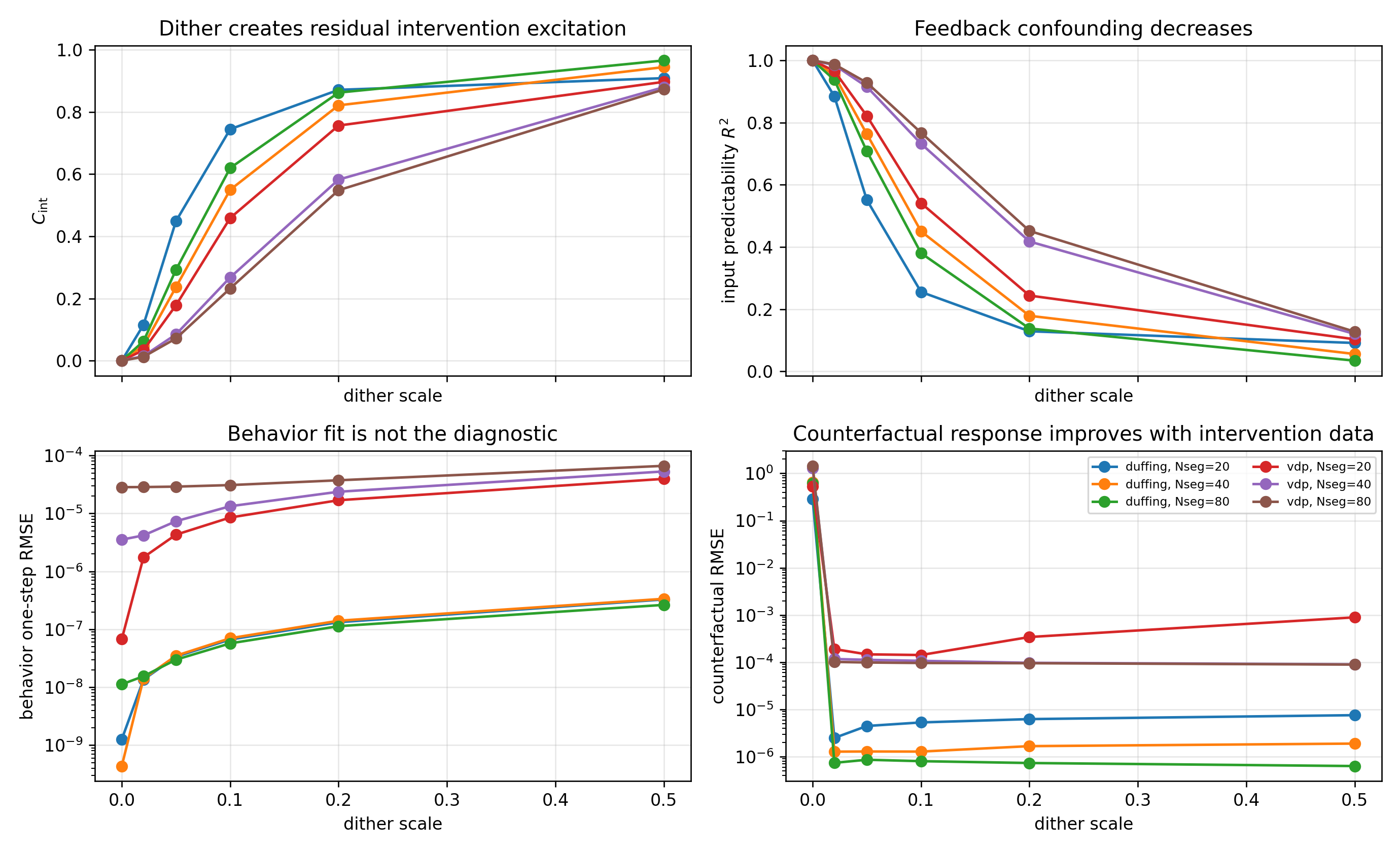}
\caption{Nonlinear feedback/dither experiments on Duffing and Van der Pol systems. Dither creates residual intervention excitation and reduces input predictability; behavior-policy one-step error remains small and is therefore not a reliable control-channel diagnostic.}
\label{fig:nonlinear}
\end{figure*}

Figure~\ref{fig:nonlinear} shows the diagnostic pattern. With deterministic feedback, the input is almost fully predictable from lifted state, $\cint$ is zero or nearly zero, and counterfactual response error is large. Adding dither sharply increases residual intervention excitation and reduces input predictability. Behavior one-step error is already small in both regimes, so behavior fit alone does not reveal the failure. Counterfactual response improves when intervention excitation is present.

\begin{table}[t]
\centering
\caption{Direct theory-validation summary. Values are medians at the largest reported budget or sample size.}
\label{tab:validation}
\begin{tabular}{llcccc}
\toprule
Experiment & condition & $\cint$ & behavior RMSE & cf. RMSE & $B$ error\\
\midrule
scalar & $\epsilon=0$, $N=800$ & $2.17{\times}10^{-35}$ & 0.0199 & 1.08 & 0.803\\
scalar & $\epsilon=0.03$, $N=800$ & $8.86{\times}10^{-4}$ & 0.0199 & 0.029 & 0.0154\\
scalar & $\epsilon=0.1$, $N=800$ & 0.00984 & 0.0198 & 0.012 & 0.0046\\
scalar & $\epsilon=0.3$, $N=800$ & 0.0885 & 0.0198 & 0.00385 & 0.00153\\
scalar & $\epsilon=1$, $N=800$ & 0.984 & 0.0198 & 0.00125 & $4.58{\times}10^{-4}$\\
Duffing & dither 0, budget 80 & 0 & $1.12{\times}10^{-8}$ & 0.610 & --\\
Duffing & dither 0.05, budget 80 & 0.292 & $2.99{\times}10^{-8}$ & $8.60{\times}10^{-7}$ & --\\
Van der Pol & dither 0, budget 80 & 0 & $2.85{\times}10^{-5}$ & 1.45 & --\\
Van der Pol & dither 0.05, budget 80 & 0.0715 & $2.91{\times}10^{-5}$ & $9.88{\times}10^{-5}$ & --\\
\bottomrule
\end{tabular}
\end{table}

\subsection{Broader acquisition grid}

The full acquisition grid is not used as the primary proof of the scaling law. It remains useful as a broader diagnostic profile. The grid spans controlled Duffing, Van der Pol, and Lorenz systems; budgets $\{8,12,20,40,80\}$; ten seeds; and ten acquisition rules including random sampling, Sobol sampling, state coverage, lifted-space and regression-space design, input-design baselines, and a certificate-aware heuristic.

\begin{figure}[t]
\centering
\includegraphics[width=0.95\columnwidth]{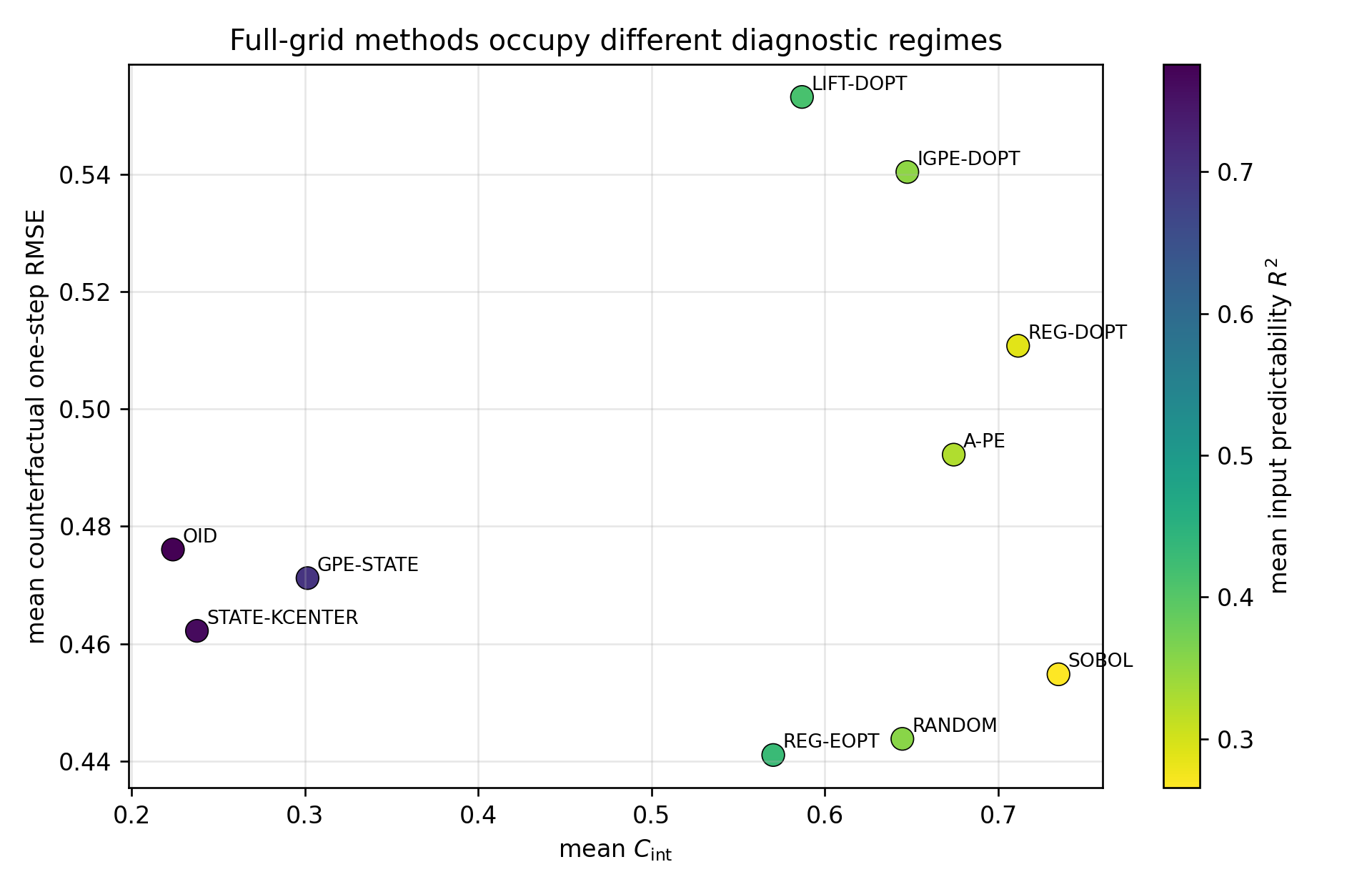}
\caption{Broad EDMDc acquisition grid. Mean $\cint$, input predictability, and counterfactual error separate different data-quality regimes; no single certificate is a universal leaderboard.}
\label{fig:v9}
\end{figure}

Figure~\ref{fig:v9} shows that acquisition rules occupy different diagnostic regimes. Some methods have high residual intervention excitation but not the lowest counterfactual error; others have low intervention excitation and high input predictability. This supports the intended interpretation: $\cint$ diagnoses a specific control-channel failure mode, while downstream errors also depend on dictionary approximation, state coverage, finite-sample noise, and nonlinear rollout effects.

\section{Practical Guidance and Limitations}

The recommended workflow is simple. Choose the dictionary, remove inactive coordinates, and standardize lifted features and inputs on the identification data. Compute $\creg$ to check whether the full regression is numerically stable. Then compute $\cint$ to check whether the input still has residual variation after conditioning on lifted state. If $\creg$ is small, the complete EDMDc fit is poorly conditioned. If $\creg$ is acceptable but $\cint$ is small, behavior prediction can still be misleading for control. The remedy is not just a longer trajectory under the same feedback law; it is data collection that breaks the lifted feedback relation.

The certificate should be compared within a fixed preprocessing and dictionary convention. Its numerical scale depends on active-coordinate selection, input dimension, standardization, and the operating region. A large value removes one confounding mechanism, but it does not validate the dictionary, guarantee long-horizon rollout, or certify nonlinear constraints. Measurement noise in the regressors would also require an errors-in-variables or instrumental-variable analysis beyond the transition-noise model used here.

\section{Conclusion}

This paper introduced control-channel informativity certificates for Koopman EDMDc under behavior-policy data. The conditional intervention certificate is the residual input covariance after projecting away active lifted-state features, equivalently the Schur complement for the EDMDc control block. It is necessary and sufficient for finite-sample sample-identifiability of the lifted control channel, and its regularized version controls closed-loop adapted-noise error bounds through self-normalized concentration. Scalar and nonlinear experiments show the central phenomenon: accurate behavior prediction can coexist with non-identifiable intervention response, while dither restores residual input information and reduces counterfactual control-channel error. The certificate is therefore best used as a targeted pre-controller diagnostic, complementary to state coverage, lifted-feature rank, and joint regression conditioning.

\appendices

\section{Proof of Proposition~\ref{prop:schur}}

Using $P_Z=Z^\top(ZZ^\top)^\dagger Z$,
\begin{equation}
\frac{1}{N}U_\perp U_\perp^\top
=\frac{1}{N}UM_ZU^\top
=\frac{1}{N}UU^\top-\frac{1}{N}UZ^\top(ZZ^\top)^\dagger ZU^\top .
\end{equation}
Since $G_{zz}^{\dagger}=N(ZZ^\top)^\dagger$, \eqref{eq:schur} follows. The rank equivalence is the Schur-complement condition for a block positive-definite matrix. The lower-right inverse block follows from the standard block-inverse formula.

\section{Proof of Theorem~\ref{thm:deterministic}}

Right-multiplying \eqref{eq:model} by $M_Z$ removes the lifted-state term:
\begin{equation}
YM_Z=B_\star UM_Z+RM_Z+EM_Z=B_\star U_\perp+RM_Z+EM_Z.
\end{equation}
If $\beta_N>0$, then $U_\perp U_\perp^\top$ is nonsingular and the residualized least-squares problem has the unique solution \eqref{eq:bhat}. This gives \eqref{eq:berr}. Conversely, if $\beta_N=0$, there exists $a\ne0$ such that $a^\top U_\perp=0$. Thus $a^\top U$ lies in the row span of $Z$, so some $h$ satisfies $a^\top U=h^\top Z$. The constructed $(A_c,B_c)$ therefore produces the same sample predictions, while \eqref{eq:counterfactual} follows by evaluation at a new pair $(z,u)$.

\section{Proof Sketch of Theorem~\ref{thm:martingale}}

The ridge normal equation gives
\begin{equation}
\widehat K_\lambda-K_\star
=E\Phi^\top V_N^{-1}+R\Phi^\top V_N^{-1}-\lambda K_\star V_N^{-1}.
\end{equation}
Under Assumption~\ref{ass:adapted}, each output coordinate of $\sum_k e_k\phi_k^\top$ is a vector martingale transform. Applying a self-normalized martingale inequality to $V_N$ bounds $\opnorm{E\Phi^\top V_N^{-1/2}}$ by the logarithmic determinant term in \eqref{eq:martbound}, up to the absolute constant $c$. The residual and ridge bias terms give the remaining numerator terms. Finally, the $B$-block of $V_N^{-1}$ has lower-right inverse governed by the regularized Schur complement $S_{u|z,\lambda}$, yielding the denominator in \eqref{eq:martbound}.

\section{Experiment Details}

The scalar experiment uses 50 seeds for every $(N,\epsilon)$ pair. Reported curves use medians across seeds. The nonlinear feedback/dither experiment uses 20 seeds per setting, trajectory segments of length 12, integration step $0.01$, and the same standardized polynomial dictionaries and ridge-selected EDMDc fitting pipeline as the acquisition grid. The counterfactual metric evaluates one-step responses under inputs not generated by the behavior policy. All result CSV files and figure-generation code are included in the project repository, while raw trajectories and large result directories are excluded from the manuscript source archive.

\section*{Data and Code Availability}

The experiment scripts, aggregated data tables, figure-generation utilities, and reproduction instructions for this manuscript are available in the project repository:
\url{https://github.com/marcowus/GPE}.
The manuscript source archive contains the LaTeX files, bibliography, compiled bibliography, README, and source figures needed to rebuild the paper; code, raw trajectories, and larger result directories are maintained as reproducibility artifacts in the repository rather than bundled into the manuscript source archive.

\bibliographystyle{IEEEtranN}
\bibliography{references}

@article{narendra1987persistent,
  title={Persistent excitation in adaptive systems},
  author={Narendra, Kumpati S and Annaswamy, Anuradha M},
  journal={International Journal of Control},
  volume={45},
  number={1},
  pages={127--160},
  year={1987},
  publisher={Taylor \& Francis},
  doi={10.1080/00207178708933715}
}

@article{willems2005note,
  title={A note on persistency of excitation},
  author={Willems, Jan C and Rapisarda, Paolo and Markovsky, Ivan and De Moor, Bart LM},
  journal={Systems \& Control Letters},
  volume={54},
  number={4},
  pages={325--329},
  year={2005},
  publisher={Elsevier},
  doi={10.1016/j.sysconle.2004.09.003}
}

@article{van2020data,
  title={Data informativity: A new perspective on data-driven analysis and control},
  author={van Waarde, Henk J and Eising, Jaap and Trentelman, Harry L and Camlibel, M Kanat},
  journal={IEEE Transactions on Automatic Control},
  volume={65},
  number={11},
  pages={4753--4768},
  year={2020},
  publisher={IEEE},
  doi={10.1109/TAC.2020.2966717}
}

@article{hjalmarsson2005experiment,
  title={From experiment design to closed-loop control},
  author={Hjalmarsson, H{\aa}kan},
  journal={Automatica},
  volume={41},
  number={3},
  pages={393--438},
  year={2005},
  publisher={Elsevier},
  doi={10.1016/j.automatica.2004.11.021}
}

@article{decock2016doptimal,
  title={D-optimal input design for nonlinear FIR-type systems: A dispersion-based approach},
  author={De Cock, Arne and Gevers, Michel and Schoukens, Johan},
  journal={Automatica},
  volume={73},
  pages={88--100},
  year={2016},
  publisher={Elsevier},
  doi={10.1016/j.automatica.2016.04.052}
}

@article{wilson2014trajectory,
  title={Trajectory synthesis for Fisher information maximization},
  author={Wilson, Alexander D and Schultz, Joshua A and Murphey, Todd D},
  journal={IEEE Transactions on Robotics},
  volume={30},
  number={6},
  pages={1358--1370},
  year={2014},
  publisher={IEEE},
  doi={10.1109/TRO.2014.2345918}
}

@article{lee2021optimal,
  title={Optimal excitation trajectories for mechanical systems identification},
  author={Lee, Taewoo and Lee, Benjamin D and Park, Frank C},
  journal={Automatica},
  volume={131},
  pages={109773},
  year={2021},
  publisher={Elsevier},
  doi={10.1016/j.automatica.2021.109773}
}

@article{kiss2024spacefilling,
  title={Space-filling input design for nonlinear state-space identification},
  author={Kiss, M{\'a}ty{\'a}s and T{\'o}th, Roland and Schoukens, Maarten},
  journal={IFAC-PapersOnLine},
  volume={58},
  number={15},
  pages={562--567},
  year={2024},
  publisher={Elsevier},
  doi={10.1016/j.ifacol.2024.08.589}
}

@article{liu2025spacefilling,
  title={On space-filling input design for nonlinear dynamic model learning: A Gaussian process approach},
  author={Liu, Yuqi and Kiss, M{\'a}ty{\'a}s and T{\'o}th, Roland and Schoukens, Maarten},
  journal={IEEE Control Systems Letters},
  volume={9},
  pages={1868--1873},
  year={2025},
  publisher={IEEE},
  doi={10.1109/LCSYS.2025.3582509}
}

@article{smits2024spacefilling,
  title={Space-filling optimized excitation signals for nonlinear system identification of dynamic processes of a diesel engine},
  author={Smits, V and Nelles, Oliver},
  journal={Control Engineering Practice},
  volume={144},
  pages={105821},
  year={2024},
  publisher={Elsevier},
  doi={10.1016/j.conengprac.2023.105821}
}

@article{williams2015datadriven,
  title={A data-driven approximation of the Koopman operator: Extending dynamic mode decomposition},
  author={Williams, Matthew O and Kevrekidis, Ioannis G and Rowley, Clarence W},
  journal={Journal of Nonlinear Science},
  volume={25},
  pages={1307--1346},
  year={2015},
  publisher={Springer},
  doi={10.1007/s00332-015-9258-5}
}

@article{proctor2016dynamic,
  title={Dynamic mode decomposition with control},
  author={Proctor, Joshua L and Brunton, Steven L and Kutz, J Nathan},
  journal={SIAM Journal on Applied Dynamical Systems},
  volume={15},
  number={1},
  pages={142--161},
  year={2016},
  publisher={SIAM},
  doi={10.1137/15M1013857}
}

@article{korda2018linear,
  title={Linear predictors for nonlinear dynamical systems: Koopman operator meets model predictive control},
  author={Korda, Milan and Mezi{\'c}, Igor},
  journal={Automatica},
  volume={93},
  pages={149--160},
  year={2018},
  publisher={Elsevier},
  doi={10.1016/j.automatica.2018.03.046}
}

@article{brunton2022modern,
  title={Modern Koopman theory for dynamical systems},
  author={Brunton, Steven L and Budi{\v{s}}i{\'c}, Marko and Kaiser, Eurika and Kutz, J Nathan},
  journal={SIAM Review},
  volume={64},
  number={2},
  pages={229--340},
  year={2022},
  publisher={SIAM},
  doi={10.1137/21M1401243}
}

@article{abraham2019active,
  title={Active learning of dynamics for data-driven control using Koopman operators},
  author={Abraham, Ian and Murphey, Todd D},
  journal={IEEE Transactions on Robotics},
  volume={35},
  number={5},
  pages={1071--1083},
  year={2019},
  publisher={IEEE},
  doi={10.1109/TRO.2019.2923880}
}

@inproceedings{boddupalli2019koopman,
  title={Koopman operators for generalized persistence of excitation conditions for nonlinear systems},
  author={Boddupalli, Nikil and Hasnain, Aqib and Nandanoori, Sai Pushpak and Yeung, Enoch},
  booktitle={Proceedings of the IEEE 58th Conference on Decision and Control (CDC)},
  pages={8106--8111},
  year={2019},
  organization={IEEE},
  doi={10.1109/CDC40024.2019.9029365}
}

@article{shang2024willems,
  title={Willems' fundamental lemma for nonlinear systems with Koopman linear embedding},
  author={Shang, Xiaoxue and Cort{\'e}s, Jorge and Zheng, Yang},
  journal={IEEE Control Systems Letters},
  year={2024},
  publisher={IEEE},
  doi={10.1109/LCSYS.2024.3522594}
}

@article{philipp2025error,
  title={Error analysis of kernel EDMD for prediction and control in the Koopman framework},
  author={Philipp, Friedrich M and Schaller, Manuel and Worthmann, Karl and Peitz, Sebastian and N{\"u}ske, Feliks},
  journal={Journal of Nonlinear Science},
  volume={35},
  pages={92},
  year={2025},
  publisher={Springer},
  doi={10.1007/s00332-025-10182-3}
}

@article{bold2025mpc,
  title={Data-driven MPC with stability guarantees using extended dynamic mode decomposition},
  author={Bold, Lea and Gr{\"u}ne, Lars and Schaller, Manuel and Worthmann, Karl},
  journal={IEEE Transactions on Automatic Control},
  volume={70},
  number={1},
  pages={534--541},
  year={2025},
  publisher={IEEE},
  doi={10.1109/TAC.2024.3431169}
}

@article{strasser2025feedback,
  title={Koopman-based feedback design with stability guarantees},
  author={Str{\"a}sser, Robin and Schaller, Manuel and Worthmann, Karl and Berberich, Julian and Allg{\"o}wer, Frank},
  journal={IEEE Transactions on Automatic Control},
  volume={70},
  number={1},
  pages={355--370},
  year={2025},
  publisher={IEEE},
  doi={10.1109/TAC.2024.3425770}
}

@article{bold2025kernel,
  title={Kernel-based Koopman approximants for control: Flexible sampling, error analysis, and stability},
  author={Bold, Lea and Philipp, Friedrich M and Schaller, Manuel and Worthmann, Karl},
  journal={SIAM Journal on Control and Optimization},
  volume={63},
  number={6},
  pages={4044--4071},
  year={2025},
  publisher={SIAM},
  doi={10.1137/24M1715945}
}

@article{rickenbach2024active,
  title={Active learning-based model predictive coverage control},
  author={Rickenbach, Remo and K{\"o}hler, Johannes and Scampicchio, Anna and Zeilinger, Melanie N and Carron, Andrea},
  journal={IEEE Transactions on Automatic Control},
  volume={69},
  number={9},
  pages={5931--5946},
  year={2024},
  publisher={IEEE},
  doi={10.1109/TAC.2024.3365569}
}

@inproceedings{abbasi2011improved,
  title={Improved algorithms for linear stochastic bandits},
  author={Abbasi-Yadkori, Yasin and P{\'a}l, D{\'a}vid and Szepesv{\'a}ri, Csaba},
  booktitle={Advances in Neural Information Processing Systems},
  volume={24},
  pages={2312--2320},
  year={2011}
}

@inproceedings{simchowitz2018learning,
  title={Learning without mixing: Towards a sharp analysis of linear system identification},
  author={Simchowitz, Max and Mania, Horia and Tu, Stephen and Jordan, Michael I. and Recht, Benjamin},
  booktitle={Proceedings of the 31st Conference on Learning Theory},
  series={Proceedings of Machine Learning Research},
  volume={75},
  pages={439--473},
  year={2018},
  publisher={PMLR}
}

\end{document}